\documentclass{article}

\usepackage{kotex}
\usepackage{latexsym,amssymb,amsmath}
\usepackage{algorithm,algorithmicx,algpseudocode}
\usepackage{calrsfs}
\usepackage{graphics}
\usepackage{graphicx}
\usepackage{setspace}
\usepackage{times}
\usepackage{ifpdf}
\usepackage{epsfig}
\usepackage{subfigure}
\usepackage[numbers]{natbib}

\parskip=\baselineskip

\providecommand{\keywords}[1]{\textbf{\textit{Keywords:}} #1}

\newtheorem{theorem}{Theorem}[section]
\newtheorem{lemma}{Lemma}[section]

\newtheorem{definition}{Definition}[section]

\newenvironment{proof}[1][Proof]{\begin{trivlist}
\item[\hskip \labelsep {\bfseries #1}]}{\end{trivlist}}

\newcommand{\qed}{\nobreak \ifvmode \relax \else
      \ifdim\lastskip<1.5em \hskip-\lastskip
      \hskip1.5em plus0em minus0.5em \fi \nobreak
      \vrule height0.75em width0.5em depth0.25em\fi}

\newcommand{\pP}{\mathbb{P}}

\begin{document}
\pagestyle{plain}
\pagenumbering{arabic}

\title{Technical Note - Exact simulation of the first passage time of Brownian motion to a symmetric linear boundary}
\author{Jong Mun Lee, Taeho Lee}
%\author{Taeho Lee, Wanmo Kang, and Jong Mun Lee}

\maketitle
%\begin{center}

%\end{center}

\begin{abstract}
We state an exact simulation scheme for the first passage time of a Brownian motion to a symmetric linear boundary. This note summarizes a part of Chapter 4 in \cite{JongMun2017} which is the unpublished doctoral dissertation of Jong Mun Lee.
\end{abstract}

\keywords{Brownian motion, hitting time simulation, exact simulation, linear boundary}

\section{Theoretical Results} 
Let $\tau_{a,b}$ be the first hitting time of Brownian motion to a symmetric and linear boundary $\pm (a+bt)$ with given nonnegative coefficients $a,b$. The following theorem from \cite{abundo02} provides an explicit expression for the distribution of $\tau_{a,b}$ including the probability $\pP[\tau_{a,b} = \infty]$.
\begin{theorem}[Abundo \cite{abundo02}] \label{thm:Abundo1}
For any $a,b>0$, the distribution of $\tau_{a,b}$ is given by
\begin{align*} \label{dist1}
\pP [ \tau_{a,b} \leq t ] &= 1 - \sum_{k \in \mathbb{Z}} (-1)^k e^{- 2 k^2 a b } \int_{(-(a+bt)+2ak)/\sqrt{t}}^{(a+bt+2ak)/\sqrt{t}} \frac{e^{-x^2 / 2}}{\sqrt{2 \pi}} dx,\\
\pP [ \tau_{a,b} < \infty ] &= 2 \sum_{k=1}^\infty (-1)^{k+1} e^{- 2 k^2 a b }.
\end{align*}
\end{theorem} 
We suggest an algorithm based on the Acceptance \& Rejection (AR) method motivated by \cite{burq08} which solves the same problem under the restriction $b=0$. We first calculate the conditional density function $f_{a,b}(t)$ under the condition $[\tau_{a,b} < \infty]$.
\begin{eqnarray}
%\begin{split}
f_{a,b} (t)  &=&  \frac{d}{dt} \pP [ \tau_{a,b} \leq t\; | \;\tau_{a,b} < \infty ] \label{eqn:desnity:tau}\\
& = & - \frac{1}{C_{a,b}} \sum_{k \in \mathbb{Z}} \frac{(-1)^k e^{- 2 k^2 a b }}{\sqrt{8 \pi t^3}}
\left((-a+bt-2ak) e^{- \frac{(a+bt+2ak)^2}{2t}} - (a-bt-2ak) e^{- \frac{(-a-bt+2ak)^2}{2t}}\right) \nonumber\\
&= & \frac{1}{C_{a,b}} \left[ \frac{a - bt}{\sqrt{2 \pi t^3}} e^{ - \frac{(a+bt)^2}{2t}} +  \sum_{k =1}^\infty \frac{(-1)^k e^{- 2 k^2 a b }}{\sqrt{2 \pi t^3}}
\left((-a+bt-2ak) e^{- \frac{(a+bt+2ak)^2}{2t}} - (a-bt-2ak) e^{- \frac{(-a-bt+2ak)^2}{2t}}\right) \right]\nonumber,
%\end{split}        
\end{eqnarray}
where $C_{a,b} =\pP[\tau_{a,b} < \infty] =  2 \sum_{k=1}^\infty (-1)^{k+1} e^{- 2 k^2 a b }$.

\noindent To apply the AR method, we need a density whose scale-up dominates the true density $f_{a,b}$. The following theorem shows that the gamma distribution can be used. 
\begin{theorem} \label{thm1}
For any $\alpha \geq \frac{1}{2}$, there exist a constant $M$ such that 
\begin{equation}\label{eq:thm1}
f_{a,b}(t) \leq M\cdot g(t ; \alpha, b^2/2)
\end{equation}
for all $t\geq0$, where $g(t:\alpha,\lambda)$ is the density function of gamma distribution given by 
\begin{equation}
g(t; \alpha,\lambda) = \frac{\lambda^\alpha t^{\alpha-1 e^{-\lambda t}}}{\Gamma(\alpha)}
\end{equation}
\begin{proof}
To prove Theorem \ref{thm1}, we analyze the tail asymptotic of $f_{a,b}(t)$. The following two lemmas analyze the left tail asymptotic and the right tail asymptotic, respectively.
\begin{lemma}[Left tail asymptotic] \label{lemma:lefttail}
For any $n \geq 1$
$$f_{a,b} (t)  = o (t^n) \text{ as $t \to 0$.}$$
\begin{proof}
Taking the positive terms and using the fact that $e^{-x} \leq m!/x^m$ for any $m \geq 0$ and $x \geq 0$, we have for $0 \leq t < a/(2b)$ and $m \geq 2$
\begin{equation*}
\begin{split}
C_{a,b} \cdot f_{a,b} (t) 
& \leq \frac{a - bt}{\sqrt{2 \pi} \sqrt{t^3}} e^{ - (a+bt)^2/(2t)} \\
& \qquad
+ \sum_{k : \text{odd}} e^{- 2 k^2 a b } \frac{1}{\sqrt{2 \pi}\sqrt{t^3}} 
(-a+bt+2ak) e^{- (-(a+bt)+2ak)^2/(2t)} \\
& \qquad
+ \sum_{k : \text{even}}  e^{- 2 k^2 a b } \frac{1}{\sqrt{2 \pi}\sqrt{t^3}} 
(a-bt+2ak) e^{- (a+bt+2ak)^2/(2t)} \\
& \leq \frac{e^{- (a-bt)^2/(2t)}}{\sqrt{2 \pi}\sqrt{t^3}} \sum_{k \geq 0 }  e^{- 2 k^2 a b } a(2k+1) \\
& \leq \frac{m! 2^m t^m}{(a-bt)^{2m}} \frac{1}{\sqrt{2 \pi}\sqrt{t^3}} \sum_{k \geq 0 }  e^{- 2 k^2 a b } a(2k+1)  = c_m t^{m-3/2} ,\\
\end{split}
\end{equation*}
where $c_m$ is some constant depends on $a,b,m$ only.
\end{proof}
\end{lemma}

\begin{lemma}[Right tail asymptotic] \label{lemma:righttail}
$$f_{a,b} (t)  = O (t^{-\frac{1}{2}} e^{-\frac{b^2}{2}t}) \text{ as $t \to \infty$.}$$
\begin{proof}
For large enough $t$ 
\begin{equation*}
\begin{split}
C_{a,b} \cdot f_{a,b} (t)  
&\leq 
\sum_{k =1}^\infty  e^{- 2 k^2 a b } \frac{2}{\sqrt{2 \pi}\sqrt{t^3}} 
\big(
(a+bt+2ak) e^{- (-(a+bt)+2ak)^2/(2t)}
\big). \\
&= \frac{2 e^{ - (a+bt)^2/(2t)} }{\sqrt{2 \pi}\sqrt{t^3}} 
\sum_{k =1}^\infty (a+bt+2ak) \exp \big( - \frac{4a^2 k^2 - 4 (a+bt) a k }{2t} - 2k^2ab \big) \\
&= \frac{2 e^{ - (a+bt)^2/(2t)} }{\sqrt{2 \pi}\sqrt{t^3}} 
\sum_{k =1}^\infty (a+bt+2ak) \exp \big( - 2ab(k^2-k) \big) \\
&= O (t^{-\frac{1}{2}} e^{-\frac{b^2}{2}t}).
\end{split}
\end{equation*}
\end{proof}
\end{lemma}
From Lemma \ref{lemma:lefttail}, we have $\lim_{t\rightarrow 0} \frac{f_{a,b}(t)}{g(t;\alpha,b^2/2)} =0$ for any $\alpha,b>0$. And from Lemma \ref{lemma:righttail}, we have $\limsup_{t\rightarrow \infty} \frac{f_{a,b}(t)}{g(t;\alpha,b^2/2)} <\infty$ for $\alpha\geq 1/2$ and thus Theorem \ref{thm1} is proved.
\end{proof}
\end{theorem}
Now we propose a prototype of our algorithm to sample $\tau_{a,b}$.
\begin{algorithm}[ht]
  \caption{Sampling of $\tau_{a,b}$( Sketch )} \label{alg:sym:proto}
  \begin{algorithmic}[1]
  	\State Generate $U \sim U(0,1)$
  	\If{ $U < C_{a,b}$} 
		\State Resolved = FALSE 
		\Repeat 
			\State Generate $V \sim \text{Gamma}(\alpha, b^2/2)$
			\If{$f_{a,b}(V) <  M \cdot g(V; \alpha, b^2/2)$}
				\State $T = V$ 
				\State Resolved = TRUE
			\EndIf
		\Until Resolved 
  	\Else 
  	\State $T= \infty$
  	\EndIf
	\State
 	\Return $T$
  \end{algorithmic}
\end{algorithm}
We call Algorithm \ref{alg:sym:proto} a prototype since the evaluation of $C_{a,b}$ and $f_{a,b}(t)$ cannot be done in finite time. Though we cannot evaluate $C_{a,b}$ and $f_{a,b}(t)$ exactly, it is possible to determine the corresponding inequalities exactly in finite time. We explore some special feature of the quantities $C_{a,b}$ and $f_{a,b}(t)$ first.
\begin{definition}[Oscillating sequence and series] 
A sequence of numbers $\{a_n\}_{n=1}^\infty$ is called an oscillating sequence if there exists a positive integer $N$ such that for $k \geq N$, 
\begin{equation*}
 -1 <  \frac{a_{k+1}-a_k}{a_k- a_{k-1}} < 0.
\end{equation*}
And we call an infinite series $\sum_{n=1}^\infty a_n$ is an oscillating series if its partial sum $s_n = \sum_{k=1}^n a_k$ is an oscillating sequence.
\end{definition}
Lemma \ref{lem:osc:prop} shows that the inequality $[S>s]$ can be  determined exactly in finite time for an oscillating series $S = \sum_{n=1}^\infty a_n$ and a constant $s\in \mathbb{R}$. The same technique is used in \cite{chen2013} which suggests an exact simulation scheme for a solution of SDE.
\begin{lemma}\label{lem:osc:prop}
Let $S= \sum_{n=1}^\infty a_n$ be an oscillating series where its oscillating property holds for $k \geq N$ and $S_n$ be its partial sum. We have $S>s$ if and only if 
\begin{equation}
S_k \wedge S_{k+1} >s  
\end{equation}
for some $k \geq N$. Conversely, We have $S<s$ if and only if 
\begin{equation}
S_k \wedge S_{k+1} <s  
\end{equation}
for some $k \geq N$. 
\end{lemma}
The following two lemmas shows that the series $C_{a,b}$ and $f_{a,b}(t)$ used in Algorithm \ref{alg:sym:proto} are oscillating series.
\begin{lemma} \label{lemma:oscil_3}
$C_{a,b} = 2 \sum_{k=1}^\infty (-1)^{k+1} e^{- 2 k^2 a b }$ is an oscillating series.
\begin{proof}
Let $a_k =(-1)^{k+1} e^{- 2 k^2 a b }$. Then we have 
$$ 0 < - \frac{a^{(n+1)} - a^{(n)}}{a^{(n)} - p^{(n-1)}} = \frac{e^{-2(n+1)^2ab}}{e^{-2n^2ab}} = e^{-2(2n+1)ab} < 1.$$
\end{proof}
\end{lemma}
\begin{lemma} \label{lemma:oscil_4}
For any $t>0$, $f_{a,b}(t)$ is an oscillating series with $N(t) = \max\left( \log{\frac{6}{4ab}}, \frac{bt}{2a},1\right) +1$
\begin{proof}
Let
\begin{equation*}
h_n (t) := e^{-2n^2ab} \big( (-a+bt+2an) e^{-(-(a+bt)+2an)^2/2t} - (a-bt+2an) e^{-(a+bt+2an)^2/2t} \big).
\end{equation*}
Then 
$$- \frac{q^{(n+1)}(t) - q^{(n)}(t) }{q^{(n)}(t) - q^{(n-1)}(t) } = \frac{h_{n+1}(t)}{h_n(t)}.$$

First, we have
\begin{equation*}
\begin{split}
h_n(t) > 0 &\Leftrightarrow (-a+bt+2an) e^{-(-(a+bt)+2an)^2/2t} > (a-bt+2an) e^{-(a+bt+2an)^2/2t}\\
&\Leftrightarrow e^{4a^2n/t + 4abn} > \frac{(a-bt+2an)}{(-a+bt+2an)} \\
&\Leftarrow e^{4abn} > \frac{(a+2an)}{(-a+2an)} .\\
\end{split}
\end{equation*}
For $n >1$, $(a+2an)/(-a+2an)$ is strictly bounded above by $3$.
Thus, $h_n(t)$ is strictly positive when $n > \max\{ \log 3 / 4ab, \; 1 \}$. 
Next, we have that for $n > \max\{ \log 3 / 4ab, \; 1 \} \vee (b/2a)t$, 
\begin{equation*}
\begin{split}
\frac{h_{n+1}(t)}{h_n(t)} < 1 &\Leftrightarrow 
e^{-2(n+1)^2ab} e^{- ( (a+bt)^2 + 4 a^2 (n+1)^2 ) / 2t } \\
&\qquad \times \big( (-a+bt+2a(n+1)) e^{4(a+bt)a(n+1)/2t} - (a-bt+2a(n+1)) e^{-4(a+bt)a(n+1)/2t} \big) \\
& \quad < 
e^{-2n^2ab} e^{- ( (a+bt)^2 + 4 a^2 n^2 ) / 2t } \big( (-a+bt+2an) e^{4(a+bt)an/2t} - (a-bt+2an) e^{-4(a+bt)an/2t} \big) \\
&\Leftrightarrow (-a+bt+2a(n+1)) e^{4(a+bt)a(n+1)/2t} - (a-bt+2a(n+1)) e^{-4(a+bt)a(n+1)/2t} \\
& \quad < (-a+bt+2an) e^{4(a+bt)a(3n+1)/2t} - (a-bt+2an) e^{4(a+bt)a(n+1)/2t} \\
&\Leftarrow 2a(2n+1) e^{4(a+bt)a(n+1)/2t} < (-a+bt+2an) e^{4(a+bt)a(3n+1)/2t} \\
&\Leftrightarrow \frac{2a(2n+1)}{(-a+bt+2an)} < e^{4(a+bt)an/t} \\
&\Leftarrow \frac{2(a+2an)}{(-a+2an)} < e^{4abn}.
\end{split}
\end{equation*}
The left-hand side of the last inequality is bounded above by $6$. And to conclude, we have that 
$$0 < - \frac{q^{(n+1)}(t) - q^{(n)}(t) }{q^{(n)}(t) - q^{(n-1)}(t) } = \frac{h_{n+1}(t)}{h_n(t)}< 1,$$
for $n \geq N(t) := \max \{ \frac{\log 6 }{4ab}, \; \frac{bt}{2a}, \; 1 \} + 1$.
\end{proof}
\end{lemma}
Since $C_{a,b},f_{a,b}(t)$ are oscillating, from lemma \ref{lem:osc:prop}, we can determine the inequalities in Algorithm \ref{alg:sym:proto} in finite time. 

\section{Algorithm}
\noindent Now we write down our algorithm to sample $\tau_{a,b}$. We divide our algorithm into two phases. In Algorithm \ref{alg:sym:step1}, we decide whether $\tau_{a,b}$ is finite or not. And in Algorithm \ref{alg:sym:step2}, we sample $\tau_{a,b}$ under the condition $[\tau_{a,b}<\infty]$. In Algorithm \ref{alg:sym:step1} and Algorithm \ref{alg:sym:step2}, we write $C^{(N)},F^{(N)}$ to denote the partial sums of $C_{a,b}$ and $f_{a,b}(t)$, respectively.
 \begin{algorithm}[ht]
  \caption{To determine whether $\tau_{a,b}$ is finite}\label{alg:sym:step1}
  \begin{algorithmic}[1]
  \State Generate $U \sim U(0,1)$
  \State $N = 1$
  \State Resolved = FALSE
  \Repeat
    \If{$C^{(N)} \vee C^{(N+1)} < U$}
      \State Resolved = TRUE
      \State Finite = FALSE
    \ElsIf{$C^{(N)} \wedge C^{(N+1)} > U$}
      \State Resolved = TRUE
      \State Finite = TRUE
    \EndIf
    \State $N = N+1$
  \Until Resolved
  \State
  \Return Finite 
  \end{algorithmic}
\end{algorithm}
\begin{algorithm}[ht]
  \caption{To generate $\tau_{a,b}$ given that $\tau_{a,b} < \infty $}\label{alg:sym:step2}
  \begin{algorithmic}[1]
 \Repeat
  \State Generate $U \sim U(0,1)$ and $V \sim \text{Gamma}(\alpha, b^2/2)$
  \State $N = \lceil N(V) \rceil$ \Comment{$N(V) = \max \{ \log 6 / 4ab, \; (b/2a) V, \; 1 \} + 1$}
  \State Resolved = FALSE
  \Repeat
    \If{$F^{(N)} (V) \vee F^{(N+1)} (V) < C_{a,b} M g(V; \alpha, b^2/2) U$}
      \State Resolved = TRUE
      \State Accept = FALSE
    \ElsIf{$F^{(N)} (V) \wedge F^{(N+1)} (V) > C_{a,b} M g(V; \alpha, b^2/2) U$}
      \State Resolved = TRUE
      \State Accept = TRUE
    \EndIf
    \State $N = N+1$
  \Until Resolved
  \Until Accept 
  \State \Return $V$
  \end{algorithmic}
\end{algorithm}

\newpage

\bibliographystyle{unsrt}
\bibliography{reference}

\begin{thebibliography}{1}

\bibitem{JongMun2017}
Jong~Mun Lee.
\newblock {\em Efficient simulation for greeks and bridges under diffusion
  models}.
\newblock PhD thesis, KAIST, 2017.

\bibitem{abundo02}
Mario Abundo.
\newblock Some conditional crossing results of brownian motion over a
  piecewise-linear boundary.
\newblock {\em Statistics \& probability letters}, 58(2):131--145, 2002.

\bibitem{burq08}
Zaeem~A Burq and Owen~D Jones.
\newblock Simulation of brownian motion at first-passage times.
\newblock {\em Mathematics and Computers in Simulation}, 77(1):64--71, 2008.

\bibitem{chen2013}
Nan Chen and Zhengyu Huang.
\newblock Localization and exact simulation of brownian motion-driven
  stochastic differential equations.
\newblock {\em Mathematics of Operations Research}, 38(3):591--616, 2013.

\end{thebibliography}

\end{document}